\documentclass[secthm,seceqn]{elsart}
\usepackage{amsfonts}
\usepackage{amssymb}
\usepackage{latexsym}

\newcommand{\defeq}{\stackrel{\rm{def}}{=}}

\journal{J.~Differential Equations}
\date{25 June 2002}

\begin{document}
\begin{frontmatter}
\title{Algebraic Solutions of the Lam\'e Equation, Revisited}
\author{Robert S. Maier\thanksref{NSF}}
\address{Depts.\ of Mathematics and Physics, University of Arizona, Tucson AZ 85721, USA}
\thanks[NSF]{Partially supported by NSF grants PHY-9800979 and
PHY-0099484.}
\begin{abstract}
A minor error in the necessary conditions for the algebraic form of the
Lam\'e equation to have a finite projective monodromy group, and hence for
it to have only algebraic solutions, is pointed~out.  [See F.~Baldassarri,
``On~algebraic solutions of Lam\'e's differential equation'', {\em
J.~Differential Equations\/} 41 (1) (1981), 44--58.]  It~is shown that if
the group is the octahedral group~$S_4$, then the degree parameter of the
equation may differ by~$\pm1/6$ from an integer; this possibility was
missed.  The omission affects a recent result on the monodromy of the
Weierstrass form of the Lam\'e equation.  [See R.~C. Churchill,
``Two-generator subgroups of ${\it SL}(2,{\bf C})$ and the hypergeometric,
Riemann, and Lam\'e equations'', {\em J.~Symbolic Computation\/} 28 (4--5)
(1999), 521--545.]  The Weierstrass form, which is a differential equation
on an elliptic curve, may have, after~all, an octahedral projective
monodromy group.
\end{abstract}
\begin{keyword}
Lam\'e equation \sep hypergeometric equation \sep projective monodromy
group \sep finite monodromy \sep algebraic solution \sep Schwarz list
\MSC 34A20 \sep 33E10 \sep 14H05
\end{keyword}
\end{frontmatter}

\section{Introduction}

The Lam\'e equation is a second-order Fuchsian differential equation.
It~may be written $L_{\ell,B}u=0$, where $L_{\ell,B}$~is the Lam\'e
operator with complex parameters $\ell$ and~$B$.  The first, the so-called
degree parameter, is often denoted~$n$, but the notation~$\ell$ is used
here, to hint at connections with Lie group representation theory.  $B$~is
an accessory parameter, which in many applications plays the role of an
eigenvalue.

The Lam\'e equation arose in a classical setting: the solution of Laplace's
equation in ellipsoidal coordinates by separation of variables.  In that
context, its solutions include the ellipsoidal harmonics.  In classical
treatments, $\ell$~is accordingly an integer, or perhaps a
half-odd-integer~\cite[Chap.~XXIII]{Whittaker27}.  The latter case arises
in a more complicated separation of variables problem (see \cite[Chap.~IX,
Ex.~4]{Poole36}, \cite[Sec.~15.1.3]{Erdelyi53}).  In~modern applications,
$\ell$~may vary continuously.  For example, the Lam\'e equation with
$\ell\in[0,2]$ has been used to compute the Hubble distance--redshift
relation in inhomogeneous, spatially flat cosmologies~\cite{Kantowski2001}.
In~that application, $\ell(\ell+1)/6\in[0,1]$ is the fraction of
inhomogeneous matter in the universe that is `dark', i.e., excluded from
observation.

Actually, several distinct equations are referred~to in the literature as
the Lam\'e equation.  We~initially consider the algebraic form, rather than
the Weierstrass or the Jacobi form.  The algebraic form is defined on the
complex projective line~${\bf P}^1({\bf C})={\bf C}\cup\{\infty\}$, with
\begin{equation}
\label{eq:algLame}
L_{\ell,B} \defeq D^2 + \frac12 \sum_{i=1}^3
\frac1{x-e_i}\,D - \frac {\ell(\ell+1)x +
B}{4\prod_{i=1}^3(x-e_i)},
\end{equation}
where $D\defeq d/dx$.  Here $\ell,B,e_1,e_2,e_3\in{\bf C}$, the $e_i$~are
distinct, and by convention, $e_1+e_2+e_3=0$.  The equation $L_{\ell,B}u=0$
has four regular singular points, three of which ($e_1,e_2,e_3$) have
characteristic exponents $0,1/2$, and one of which~($\infty$) has exponents
$-\ell/2$, $(\ell+1)/2$.  So the algebraic-form Lam\'e equation is a
special case of the Heun equation, which is the general second-order
Fuchsian equation on~${\bf P}^1({\bf C})$ with four singular points.

Via the map $(x,y)\mapsto x$, the line ${\bf P}^1({\bf C})$~is doubly
covered by the elliptic curve $y^2=4x^3-g_2x-g_3$, where the invariants
$g_2,g_3\in\bf C$, at~least one of which is nonzero, are defined by
$4x^3-g_2x-g_3\equiv 4\prod_{i=1}^3(x-e_i)$.  This curve will be denoted
$E_{g_2,g_3}$.  $L_{\ell,B}$~can be pulled back to a differential operator
$L_{\ell,B,g_2,g_3}$ that acts on~$E_{g_2,g_3}$.  The pullback has only one
singular point, namely $(x,y)=(\infty,\infty)$, which from a geometric
point of view is why the Lam\'e case of the Heun equation is important.
The pulled-back equation ${L_{\ell,B,g_2,g_3}u=0}$ on~$E_{g_2,g_3}$ is the
Weierstrass form, which is discussed in Section~\ref{sec:corrections}.
Indirect references to the elliptic curve interpretation occur elsewhere,
since when studying~$L_{\ell,B}$, we classify various situations by
supplying the corresponding value of~$J$, Klein's modular function (also
known as Klein's absolute invariant).  Recall that $J\defeq
g_2^3/\Delta\in{\bf C}$, where $\Delta\defeq g_2^3-27g_3^2\neq0$ is the
modular discriminant.  Iff two elliptic curves specified by $g_2,g_3$ have
the same value of~$J$, they are birationally equivalent, e.g., homeomorphic
as complex manifolds~\cite[Sec.~5.3]{Cohn85}.

The determination of all quadruples $\ell,B,g_2,g_3$ for which $L_{\ell,B}$
has only algebraic functions in its kernel is an unsolved problem.  The
nonclassical case ${2\ell\notin\bf Z}$ of this problem is perhaps the most
tractable.  Singer~\cite{Singer93} and Morales-Ruiz and
Sim\'o~\cite[Lemma~1]{Morales96} mention an unpublished result of Dwork
that for any fixed~$\ell$ for which $2\ell\notin\bf Z$, if $e_1$~is fixed,
then there are only a finite number of pairs $e_3,B$ for which all
solutions of $L_{\ell,B}u=0$ are algebraic.  In~essence, for each
$\ell\notin(1/2){\bf Z}$ there are only a finite number of `algebraic'
pairs $J,B$; though for this statement to make sense, $B$~would need to be
redefined in a scale-invariant way, constant on each elliptic curve
isomorphism class.

The difficulty of finding all $\ell,B,g_2,g_3$ for which the Lam\'e
equation $L_{\ell,B}u=0$ has only algebraic solutions contrasts with the
classical solution of the corresponding problem for the hypergeometric
equation ${L_{\lambda,\mu,\nu}v=0}$, the canonical second-order Fuchsian
equation on ${\bf P}^1({\bf C})$ with three singular points.  Here
\begin{equation}
L_{\lambda,\mu,\nu}
\defeq \frac{d^2}{dz^2}+
\frac{1-\lambda^2}{4z^2} + \frac{1-\mu^2}{4(z-1)^2} 
+\frac{\lambda^2 + \mu^2 - 1 -\nu^2}{4z(z-1)}
\end{equation}
is the (normal-form) hypergeometric operator with exponent differences
$\lambda,\mu,\nu\in{\bf C}$, and the singular points on~${\bf P}^1({\bf
C})$ (coordinatized by~$z$) are $z=0,1,\infty$.  It~is a classical result
of Schwarz that if $\lambda,\mu,\nu\notin\bf Z$, then
$L_{\lambda,\mu,\nu}v=0$ will have only algebraic solutions iff a suitably
normalized version of $\lambda,\mu,\nu$ (regarded as an unordered triple)
appears on a certain list.  This is the famous `Schwarz list', which has 15
entries, numbered I--XV\null. (See~\cite[Sec.~30]{Poole36},
\cite[Sec.~2.7.2]{Erdelyi53},~\cite{Osgood97}.)  The case when one of
$\lambda,\mu,\nu$ is an integer is degenerate, and can be handled by other
means (it~has its own list).  To~each list entry there corresponds a finite
group, to which the projective monodromy group $G(L_{\lambda,\mu,\nu})$,
which will be a finite subgroup of the M\"obius group ${\it PGL}(2,{\bf
C})$, is necessarily isomorphic.  The possible groups are cyclic~($C_n$,
$n\ge1$), dihedral~($D_{n}$, $n\ge2$), tetrahedral~($A_4$),
octahedral~($S_4$), and icosahedral~($A_5$).

Klein's theory of pullbacks of Fuchsian operators grew out of Schwarz's
classification theory.  Associated to any second-order Fuchsian
operator~$F$ on an algebraic curve over~${\bf C}$ is a projective monodromy
group~$G(F)\leq{\it PGL}(2,{\bf C})$.  Klein showed that $G(F)$~will be
finite, which is almost enough to ensure that $Fu=0$~has only algebraic
solutions, iff $F$~is a (weak) pullback from ${\bf P}^1({\bf C})$ of some
$L_{\lambda,\mu,\nu}$, where $\lambda,\mu,\nu$ belongs to a small sublist,
called the `basic Schwarz list'.  (Other list entries can be omitted since
they are redundant: they themselves correspond to pullbacks.)  Necessarily
$G(F)\leq G(L_{\lambda,\mu,\nu})$; and in~fact, there is at~least one
$\lambda,\mu,\nu$ on the basic Schwarz list, with corresponding pullback,
such that $G(F)=G(L_{\lambda,\mu,\nu})$.  If the pullback is known
explicitly, $G(F)$ may readily be computed, and the solutions of $Fu=0$ may
be computed too.  All solutions will be algebraic, provided the Wronskian
of~$F$ is algebraic.  The proofs of Klein were modernized by Baldassarri
and Dwork in~\cite{Baldassarri79,Baldassarri80}.

In a remarkable paper, Baldassarri~\cite{Baldassarri81} applied Klein's
theory to the Lam\'e equation.  By determining necessary conditions for the
existence of a pullback of $L_{\ell,B}$ from each possible
$L_{\lambda,\mu,\nu}$, Baldassarri derived a necessary condition for
$L_{\ell,B}u=0$ to have only algebraic solutions, and also necessary
conditions for $G(L_{\ell,B})$ to be each possible finite subgroup of the
M\"obius group.  It~is a classical result that $G(L_{\ell,B})$ is never
cyclic, and can be dihedral only if $2\ell\in\bf Z$.  Moreover, in the
nonclassical case $2\ell\notin\bf Z$, it cannot be dihedral.  Baldassarri
showed that if $2\ell\notin\bf Z$, all solutions of $L_{\ell,B}u=0$ can be
algebraic only if one of $\ell\pm1/10$, $\ell\pm1/6$, $\ell\pm1/4$,
or~$\ell\pm3/10$ is an integer.  Moreover, $G(L_{\ell,B})$ cannot be
tetrahedral, so if it is finite, it must be octahedral or icosahedral.

Unfortunately, \cite{Baldassarri81} errs in its treatment of the octahedral
case.  In Theorem~\ref{thm:main}, we restate the conditions
of~\cite{Baldassarri81} with the following correction: For $G(L_{\ell,B})$
to be octahedral, it is necessary that one of $\ell\pm1/6$ or $\ell\pm1/4$
be an integer, but not that one of $\ell\pm1/4$ be an integer.
We~discovered the need for this correction while examining the implications
for Lam\'e monodromy of~\cite{Maier03}, which in~effect classifies all {\em
strong\/} pullbacks of the hypergeometric to the Heun equation.  Pulling
back `algebraic' $L_{\lambda,\mu,\nu}$ via the quadratic and cubic cyclic
maps treated in~\cite{Maier03} yields useful examples of Lam\'e operators
with only algebraic functions in their kernels, including a counterexample
to the necessary condition of~\cite{Baldassarri81}.  The counterexample
appears in Proposition~\ref{prop:main}, and explicit formul\ae\ for the
solutions of a number of interesting Lam\'e equations with projectively
finite monodromy are given in Section~\ref{sec:explicit}.

The corrected necessary condition for $G(L_{\ell,B})$ to be octahedral
overlaps with the necessary condition that it be icosahedral, which is that
one of $\ell\pm1/10$, $\ell\pm1/6$, or $\ell\pm3/10$ be an integer.  For
example, $\ell=1/6$ is both an octahedral and an icosahedral alternative.
It~follows from Propositions \ref{prop:main} and~\ref{prop:partial} that
there are Lam\'e operators with $\ell=1/6$ of both the octahedral and
icosahedral types.  This implies that in the nonclassical $2\ell\notin\bf
Z$ case, finite projective monodromy is not determined uniquely by~$\ell$.

Churchill~\cite{Churchill99} studied the monodromy of the Weierstrass-form
Lam\'e equation $L_{\ell,B,g_2,g_3}u=0$ on the elliptic
curve~$E_{g_2,g_3}$, and employed the results of~\cite{Baldassarri81} to
derive similar results on the projective monodromy group
$G(L_{\ell,B,g_2,g_3})$.  In~particular, he deduced that it cannot be
octahedral.  Unfortunately, this deduction is invalidated by the error
in~\cite{Baldassarri81} and the consequent nonuniqueness.  In
Section~\ref{sec:corrections}, we provide details, including
Theorem~\ref{thm:mainmain}, a corrected theorem on $G(L_{\ell,B,g_2,g_3})$
and its relation to $G(L_{\ell,B})$.  We~also give an example of an
equation $L_{\ell,B,g_2,g_3}u=0$ with octahedral projective monodromy.

\section{Preliminaries}
\label{sec:prelims}

The following definitions and results are fairly
standard~\cite{Baldassarri79,Baldassarri80}, but are included to make this
paper self-contained.  Suppose $C$ is a nonsingular algebraic curve
over~$\bf C$ with function field~$K/\bf C$, and that $D$~is a nontrivial
derivation of~$K/{\bf C}$.  (For example, $C={\bf P}^1({\bf C})$, with
$K={\bf C}(x)$, the field of rational functions, and $D$~the usual
derivation $d/dx$.)  Consider the monic second-order operator
\begin{equation}
\label{eq:monic}
L= D^2 + {\mathcal A}\cdot D + {\mathcal B}
\end{equation}
where ${\mathcal A},{\mathcal B}\in K$.  Let $\{P_1,\ldots,P_r\}$ be its
set of singular points, which comprises the poles of ${\mathcal A}$
and~${\mathcal B}$, and possibly the point at infinity; and let $P$~be an
ordinary point.  A~${\it GL}(2,{\bf C})$ monodromy representation of the
fundamental group of the punctured curve,
$\pi_1(C\setminus\{P_1,\ldots,P_r\};P)$, is obtained by analytically
continuing any two linearly independent function elements $u_1,u_2$ around
closed paths that issue from~$P$\null.  Its image in ${\it GL}(2,{\bf C})$
is the monodromy group of~$L$ (its~isomorphism class is independent of the
choice of $u_1,u_2$ and~$P$)\null.  The image of the monodromy group in
${\it PGL}(2,{\bf C})$, obtained by quotienting~out its intersection with
${\bf C}\setminus\{0\}$, is the projective monodromy group $G(L)$, the
group of monodromies of the ratio~$u_2/u_1$.

Iff $G(L)$ is finite, any ratio of independent solutions of $Lu=0$ will be
algebraic over~$K$, with Galois group~$G(L)$.  Let $\tau$ be such a ratio.
By calculation, if ${\mathcal A}=0$, then $u_1\defeq1/\sqrt{D\tau}$ and
$u_2\defeq\tau/\sqrt{D\tau}$ will satisfy $Lu_i=0$.  Moreover, these
$u_1,u_2$ are independent.  So if ${\mathcal A}=0$, all solutions of $Lu=0$
are algebraic over~$K$ iff $G(L)$~is finite.  ${\mathcal A}=0$ can be
weakened to the condition that the Wronskian $w=w(L)$, defined locally
on~$C$ by $Dw+{\mathcal A}\cdot w=0$, be algebraic over~$K$.  This is
because
\begin{eqnarray}
\label{eq:similarity}
\hat L &=& 1/\sqrt{w} \circ L \circ \sqrt{w}\\
&=& D^2 - D{\mathcal A}/2 - {\mathcal A}^2/4 + {\mathcal B}\nonumber
\end{eqnarray}
equals $D^2+\hat{\mathcal A}\cdot D + \hat{\mathcal B}$ with $\hat{\mathcal
A}=0$, i.e., is of `normal form'.  The groups $G(L)$ and~$G(\hat L)$ are
isomorphic, and $Lu=0$ iff $\hat L(u/\sqrt w)=0$.  That~is, the solution
space of $Lu=0$ is spanned by $\sqrt{w}/\sqrt{D\hat\tau}$ and
$\sqrt{w}\hat\tau/\sqrt{D\hat\tau}$, where $\hat\tau$~is any ratio of
solutions of $\hat L\hat u=0$; $\hat\tau$~is algebraic iff $G(L)$~is
finite.  So if $w(L)$~is algebraic, $Lu=0$ has a full set of algebraic
solutions iff $G(L)$ is finite.

Let $\xi:C\to C'$ be a rational map of algebraic curves, where $C'$ is
another nonsingular algebraic curve over~$\bf C$, with its own function
field~$K'/\bf C$ and nontrivial derivation~$D'$.  If $L$~is as
in~(\ref{eq:monic}), and $L'=(D')^2+{\mathcal A}'\cdot D'+{\mathcal B}'$, with
${\mathcal A}',{\mathcal B}'\in K'$, is a similar monic second-order operator
on~$C'$, then $L$ is said to be a strong pullback of~$L'$ if there are
independent solutions $u_1,u_2$ and $u_1',u_2'$ of $L,L'$ respectively,
such that $u_i=u_i'\circ\xi$.  For example, if $x$~is the coordinate on~$C$
and $C'={\bf P}^1({\bf C})$ is coordinatized by~$z$, so that $z=\xi(x)$~is
a rational function on~$C$, and $L'=D_z^2+{\mathcal B}'$, then the strong
pullback of~$L'$ is
\begin{eqnarray}
\label{eq:strongpullback}
&&(d\xi/dx)^2 \left[ d^2/d\xi^2 + {\mathcal B}'(\xi)\right] =\\
&&\qquad\frac{d^2}{dx^2} - \frac{d^2\xi/dx^2}{d\xi/dx}\frac{d}{dx} +
(d\xi/dx)^{2}{\mathcal B}'(\xi(x)),\nonumber
\end{eqnarray}
where the prefactor $(d\xi/dx)^2$ ensures monicity.

If $L,M$ are monic second-order operators on~$C$ (resp.~$C'$), $L$~is said
to be projectively equivalent to~$M$ (written $L\sim M$) iff any ratio of
independent solutions of~$Mu=0$ is a ratio of solutions of~$Lu=0$, i.e.,
iff $M=h^{-1}\circ L\circ h$ for some $h\in K$ (resp.~$K'$).  Note that if
$M$~is normal-form, then $h\propto\sqrt{w(L)}$ as in~(\ref{eq:similarity}),
and $M$~is uniquely determined by~$L$.

If $L,L'$ are monic second-order operators on $C,C'$, $L$~is said to be a
{\em weak\/} pullback of~$L'$ (hereafter, a pullback) if there are $\tau\in
K$, $\tau'\in K'$, ratios of independent solutions of $Lu=0$, $L'u'=0$,
with $\tau=\tau'\circ\xi$.  That~is, $L\sim M$ and $L'\sim M'$, with $M$~a
strong pullback of~$M'$.  Pullbacks are not unique, though there is a
unique normal-form pullback.

\begin{lem}
\label{lem:firstlem}
$L=D^2 + {\mathcal A}\cdot D + {\mathcal B}$ on~$C$ is a pullback by
$\xi:C\to{\bf P}^1({\bf C})$ of the normal-form operator $L'=D_z^2 +
{\mathcal B}'$ on ${\bf P}^1({\bf C})$ iff
\begin{eqnarray}
\label{eq:inlemma}
&&- D{\mathcal A}/2 -{\mathcal A}^2/4 + {\mathcal B}=\\
&&\qquad\frac12\frac{d}{dx}\left(\frac{d^2\xi/dx^2}{d\xi/dx}\right)
-\frac14\left(\frac{d^2\xi/dx^2}{d\xi/dx}\right)^2 + (d\xi/dx)^{2}B'(\xi(x)).
\nonumber
\end{eqnarray}
If this is the case, the solution space of $Lu=0$ is spanned by
\begin{equation}
\label{eq:inlemma2}
\frac{\sqrt{w(L)}}{\sqrt{D(\tau'\circ\xi)}},
\qquad
\frac{\sqrt{w(L)}(\tau'\circ\xi)}{\sqrt{D(\tau'\circ\xi)}},
\end{equation}
where $\tau'$ is any ratio of independent solutions of $L'u'=0$.
\end{lem}
\begin{pf}
The strong pullback of~$L'$ is given by~(\ref{eq:strongpullback}), and
according to the formula~(\ref{eq:similarity}), the unique normal-form weak
pullback of~$L'$ will be $D^2+\hat {\mathcal B}$, where $\hat {\mathcal
B}$~is defined as the right-hand side of~(\ref{eq:inlemma}).  But as
computed in~(\ref{eq:similarity}), an operator of the form $D^2+\hat
{\mathcal B}$ is projectively equivalent to~$L$ iff $\hat {\mathcal
B}=-D{\mathcal A}/2-{\mathcal A}^2/4+{\mathcal B}$.  The final statement
follows from the above remarks about the solution space of $Lu=0$ in
relation to that of $(D^2+\hat {\mathcal B})\hat u=0$. \qed
\end{pf}

We now specialize to operators $F=D^2+{\mathcal A}\cdot D+{\mathcal B}$
on~$C$ that are Fuchsian, i.e., have two characteristic exponents
$\alpha_{i,1},\alpha_{i,2}\in\bf C$ (which may be the same) at each
singular point~$P_i$.  If $\alpha_{i,1}\neq\alpha_{i,2}$, this means $Fy=0$
has solutions $y_{i,j}$, $j=1,2$, at~$P_i$ that are of the form
$t^{\alpha_{i,j}}$~times an invertible function of~$t$, where $t$~is a
local uniformizing parameter (if $\alpha_{i,1}=\alpha_{i,2}$, one solution
may be logarithmic).  The exponent differences
$\rho(F,P_i)=\alpha_{i,1}-\alpha_{i,2}$ are defined up~to sign; when
$\rho\in\bf R$, the convention $\rho\ge0$ will be adhered~to.  At~each
ordinary point, the exponents are $0,1$, so the exponent difference is
unity.

Theorem~\ref{thm:pullback} is Klein's pullback theorem, taken
from~\cite[Thm.~1.8]{Baldassarri80}.  The auxiliary Table~\ref{tab:list} is
the basic Schwarz list of exponent differences $\lambda,\mu,\nu$.  The
hypergeometric operator $L_{\lambda,\mu,\nu}$ corresponding to each row has
a full set of algebraic solutions, and there is a ratio~$\tau'$ of
solutions which, as an algebraic function from $z\in{\bf P}^1({\bf C})$ to
$\tau'\in{\bf P}^1({\bf C})$, is the inverse of a single-valued, i.e.,
rational, `polyhedral function' $z=z(\tau')$.  These are tabulated in the
final column, which is adapted from \cite[Sec.~31]{Poole36}
and~\cite[Sec.~14.3]{Sansone69}.  Each is automorphic under the
corresponding finite subgroup of the M\"obius group.

\begin{table}
\caption{The basic Schwarz list.}
\label{tab:list}
\begin{center}                                           
\begin{tabular}{c|l|c|l} 
\hline
Case & {\hfil$\lambda,\mu,\nu$\hfil} & {\hfil Group\hfil} & {\hfil Solution
ratio inverse, $z=z(w)$\hfil}\\
\hline
--- & $1/n,1,1/n$ & $C_n$ & $w^n$\\
I & $1/2,1/2,1/n$ & $D_n$ & ${\displaystyle\frac{(w^n+1)^2}{4w^n}}$\\
II & $1/2,1/3,1/3$ & $A_4$ & ${\displaystyle\frac{12\sqrt{-3}\,w^2(w^4-1)^2}{(w^4+2\sqrt{-3}\,w^2+1)^3}}$\\
IV & $1/2,1/3,1/4$ & $S_4$ & ${\displaystyle\frac{-(w^{12}-33w^8-33w^4+1)^2}{108\,w^4(w^4-1)^4}}$\\
VI & $1/2,1/3,1/5$ & $A_5$ & 
${\displaystyle\frac{[w^{30}+522(w^{25}-w^5)-10005(w^{20}+w^{10})+1]^2}{1728\,w^5(w^{10}+11w^5-1)^5}}$\\
\hline
\end{tabular}                                                        
\end{center}                                                       
\end{table}

\begin{thm}
\label{thm:pullback}
Let $F=D^2+{\mathcal A}\cdot D+{\mathcal B}$ be a Fuchsian operator on~$C$,
with ${\mathcal A},{\mathcal B}\in K$, and suppose that $G(F)$ is finite.
There is a unique $\lambda,\mu,\nu$ on the basic Schwarz list such that
$G(F)$ is isomorphic to $G(L_{\lambda,\mu,\nu})$ and $F$~is a pullback of
$L_{\lambda,\mu,\nu}$ by some rational map $\xi:C\to{\bf P}^1({\bf C})$,
where $\xi$~is unramified over ${\bf P}^1({\bf C})\setminus\{0,1,\infty\}$.
Moreover, if for any $\lambda,\mu,\nu$ on the list, $F$~is a pullback of
$L_{\lambda,\mu,\nu}$, then $G(F)$ is isomorphic to a subgroup of
$G(L_{\lambda,\mu,\nu})$.
\end{thm}

If $G(F)$ is finite and the Wronskian $w(F)$ is algebraic, $Fu=0$ will have
a full set of algebraic solutions; and if $\lambda,\mu,\nu$ and the
pullback map $\xi:C\to{\bf P}^1({\bf C})$, which are guaranteed to exist by
Theorem~\ref{thm:pullback}, are known, a basis for the solution space of
$Fu=0$ may be computed from~(\ref{eq:inlemma2}), in which $\tau'=\tau'(z)$
is the inverse of the polyhedral function in the final column of the table.

It is worth noting that as algebraic functions, the possible~$\tau'$ are
quite special.  Each is ramified over $z=0,1,\infty$, at~most, and the
ramification order of each of the points in $(\tau')^{-1}(0)$,
$(\tau')^{-1}(1)$, $(\tau')^{-1}(\infty)$ is $1/\lambda-1$, $1/\mu-1$,
$1/\nu-1$, respectively.  That~is, if $m$~denotes the mapping degree of
$z=z(\tau')$, i.e., $m=|G(L_{\lambda,\mu,\nu})|$, the projective monodromy
of~$L_{\lambda,\mu,\nu}v=0$ around each of the singular points
$z=0,1,\infty$ is always a restricted sort of permutation of the
$m$~branches of~$\tau'$, the cycle decomposition of which comprises,
respectively, $\lambda m$~cycles of length~$1/\lambda$, $\mu m$~cycles of
length~$1/\mu$, and $\nu m$~cycles of length~$1/\nu$.

\begin{lem}
\label{lem:exptoexp}
Suppose the Fuchsian operator $F=D^2+{\mathcal A}\cdot D+{\mathcal B}$ is a
pullback of~$L_{\lambda,\mu,\nu}$ via $\xi:C\to{\bf P}^1({\bf C})$.  The
exponent difference $\rho(F,P)$ at any $P\in C$ equals $h$~times the
exponent difference $\rho(L_{\lambda,\mu,\nu},\xi(P))$, if $h$~is the
multiplicity with which $P$~is mapped to $\xi(P)$, i.e., $1$~plus the
ramification order of $\xi$ at~$P$.
\end{lem}
\begin{pf}
Consider the series expansions of solution ratios $\tau,\tau'$ of $Fy=0$,
$L_{\lambda,\mu,\nu}y=0$ at $P,\xi(P)$, respectively.  Each is of the
form~$t^{\rho}$ times an invertible function of~$t$, where $t$~is a
uniformizing parameter; and locally, $\xi(t)\sim t^h$. \qed
\end{pf}

Lemma~\ref{lem:exptoexp} constrains the Fuchsian operators~$F$ to which
Theorem~\ref{thm:pullback} can be applied, i.e., $F$~for which $G(F)$ is
finite.  For example, there must be a row of Table~\ref{tab:list} such that
each of the singular point exponent differences $\{\rho(F,P_i)\}_{i=1}^r$
is an integer multiple of one of the corresponding $\lambda,\mu,\nu$.
It~also constrains the monodromy at each~$P_i$.  Suppose WLOG that
$\xi(P_i)=0$.  The projective monodromy of $L_{\lambda,\mu,\nu}v=0$ around
$z=0$ permutes the $m$~branches of~$\tau'$, and the cycle decomposition of
the permutation comprises $\lambda m$~cycles of length~$1/\lambda$.  So the
projective monodromy of $Fu=0$ around~$P_i$ must be isomorphic to an
integer power of such a permutation.  Together with the fact that $G(F)$,
the group of permutations of the branches of~$\tau$ which is generated by
these monodromies, must be identical to the Galois group of~$\tau$ over~$K$
(rather than being a proper subset of~it), this imposes substantial
constraints.

The following lemma will be used in the next section.

\begin{lem}
\label{lem:lookahead}
If $L_{\ell,B}$ is an algebraic-form Lam\'e operator with finite projective
monodromy group, so that it is a pullback of some $L_{\lambda,\mu,\nu}$ on
the basic Schwarz list by a rational map $\xi:{\bf P}^1({\bf C})\to{\bf
P}^1({\bf C})$ of the sort guaranteed to exist by
Theorem~\ref{thm:pullback}, then provided $\ell+1/2\notin\bf Z$, $\xi$~must
map the set of singular points $\{e_1,e_2,e_3,\infty\}$ into
$\{0,1,\infty\}$.
\end{lem}

\begin{pf}
The only ramification points of~$\xi$ are above $z=0,1,\infty$.  So if
$\xi(P)\notin\{0,1,\infty\}$,
$\rho(L_{\ell,B},P)=\rho\left(L_{\lambda,\mu,\nu},\xi(P)\right)=1$ by
Lemma~\ref{lem:exptoexp}.  Since $\rho(L_{\ell,B},e_i)=1/2$ and
$\rho(L_{\ell,B},\infty)=\pm(\ell+1/2)$, the claim follows. \qed
\end{pf}

\section{Key Results}
\label{sec:result}

\begin{thm}
\label{thm:main}
The equation $L_{\ell,B}u=0$ on~${\bf P}^1({\bf C})$ has a full set of
algebraic solutions, i.e., solutions algebraic over~${\bf C}(x)$, iff
$G(L_{\ell,B})$ is finite.  In the nonclassical case $2\ell\notin\bf Z$,
$G(L_{\ell,B})$ is finite iff it is octahedral {\rm(}i.e., isomorphic
to~$S_4${\rm)}, in which case $\ell$~must equal $n\pm1/6$ or~$n\pm1/4$,
with $n$~an integer; or icosahedral {\rm(}i.e., isomorphic to~$A_5${\rm)},
in which case $\ell$~must equal $n\pm1/10$, $n\pm1/6$, or $n\pm3/10$, with
$n$ an integer.
\end{thm}

\begin{pf}
The Wronskian $w(L_{\ell,B})$ equals $\prod_{i=1}^3 (x-e_i)^{-1/2}$, which
is algebraic; so $L_{\ell,B}u=0$ having a full set of algebraic solutions
is equivalent to finiteness of~$G(L_{\ell,B})$.  The necessary conditions
on~$\ell$ come from conditions imposed by Lemma \ref{lem:exptoexp} on
pullbacks of~$L_{\ell,B}$ from $L_{\lambda,\mu,\nu}$ on the basic Schwarz
list, since such a pullback is guaranteed to exist by
Theorem~\ref{thm:pullback}.  As~the final sentence of that theorem
acknowledges, a Fuchsian operator~$F$ can be a pullback
of~$L_{\lambda,\mu,\nu}$ with $G(F)$ isomorphic to a proper subgroup
of~$G(L_{\lambda,\mu,\nu})$, rather to $G(L_{\lambda,\mu,\nu})$.
To~compensate, one must consider the various $L_{\lambda,\mu,\nu}$
`in~order'.  The rows of Table~\ref{tab:list} are ordered so that if
$G_2$~appears in a later row than~$G_1$, then $G_2$~is not isomorphic to a
subgroup of~$G_1$.

The analysis begins with the tetrahedral row, since it is a classical
result that if $2\ell\notin\bf Z$, $G(L_{\ell,B})$ cannot be cyclic or
dihedral.  If $G(L_{\ell,B})$ is tetrahedral, $L_{\ell,B}$~must be a
pullback of $L_{1/2,1/3,1/3}$.  Since $L_{\ell,B}$~has exponent differences
$1/2$,~$1/2$, $1/2$, $\pm(\ell+1/2)$ at $x=e_1,e_2,e_3,\infty$,
respectively, it~follows from Lemma~\ref{lem:lookahead} and
Lemma~\ref{lem:exptoexp} that $\xi$~must map $e_1,e_2,e_3$ to~$0$, and
$\infty$~to~$1$, resp.~$\infty$.  Hence $\xi^{-1}(\infty)$,
resp.~$\xi^{-1}(1)$, must comprise only ordinary points with exponent
differences equal to unity.  By Lemma~\ref{lem:exptoexp}, each point
in~$\xi^{-1}(\infty)$, resp.~$\xi^{-1}(1)$, must be mapped triply
to~$\infty$.  So $3\mid\deg\xi$.  This can be combined with the prediction
of the `degree formula' of Baldassarri and
Dwork~\cite[Lemma~1.5]{Baldassarri79}, which is derived from the Hurwitz
genus formula.  If~$F$, a second-order Fuchsian operator on an algebraic
curve~$C$ with genus~$g$, has exponent differences~$\{\rho_i\}$ and is a
pullback by a rational function~$\xi$ from~$F'$, a Fuchsian operator
on~${\bf P}^1({\bf C})$ with exponent differences~$\{\rho'_i\}$, then
\begin{equation}
\label{eq:degreeformula}
\biggl[2-2g+\sum_i(\rho_i-1)\biggr] =
(\deg\xi)\biggl[2+\sum_i(\rho'_i-1)\biggr].
\end{equation}
The degree formula~(\ref{eq:degreeformula}) yields
$\pm(\ell+1/2)-1/2=(\deg\xi)/6$ when applied to $F=L_{\ell,B}$,
$F'=L_{1/2,1/3,1/3}$.  In~conjunction with $3\mid\deg\xi$, this contradicts
$2\ell\notin\bf Z$.  [This ruling~out of the tetrahedral alternative is
taken from~\cite[Prop.~3.1]{Baldassarri81}.]

If $G(L_{\ell,B})$ is octahedral, $L_{\ell,B}$~must be a pullback of
$L_{1/2,1/3,1/4}$.  The point $x=\infty$ cannot be mapped to~$0$, since by
Lemma~\ref{lem:exptoexp} that would imply that $\ell+1/2$ is an integer
multiple of~$1/2$, which is a contradiction.  However, it~can be mapped
to~$1$, in which case $\ell+1/2$ must be an integer multiple of~$1/3$, or
to~$\infty$, in which case $\ell+1/2$ must be an integer multiple of~$1/4$.
That~is, $\ell$~must equal $n\pm1/6$ or $n\pm1/4$, with $n$~an integer.
[The possibility that $\xi(\infty)=1$ was erroneously ruled~out
in~\cite[Sec.~3]{Baldassarri81}, by an argument based on the incorrect
assumption that $\xi(e_i)$ must equal~$0$ for all~$i$.]

If $G(L_{\ell,B})$ is icosahedral, $L_{\ell,B}$~must be a pullback of
$L_{1/2,1/3,1/5}$.  As in the octahedral case, $x=\infty$ cannot be mapped
to~$0$.  It~can be mapped to~$1$, in~which case $\ell+1/2$ must be an
integer multiple of~$1/3$, or to~$\infty$, in~which case $\ell+1/2$ must be
an integer multiple of~$1/5$.  That~is, $\ell$~must equal $n\pm1/6$, with
$n$~an integer, or $n\pm1/10$ or~$n\pm3/10$, with $n$~an integer. \qed
\end{pf}

According to Propositions \ref{prop:main} and~\ref{prop:partial} below, the
five alternatives listed in Theorem~\ref{thm:main} can each be realized.

\begin{defn}
The harmonic case is the case when $J=1$, i.e., when $g_3=0$, so that the
unordered set $\{e_1,e_2,e_3\}$ comprises three equally spaced collinear
points in~$\bf C$, i.e., is of the form $\alpha\{-1,0,1\}$.  The
equianharmonic case is the case when $J=0$, i.e., when $g_2=0$, so that
$\{e_1,e_2,e_3\}$ is the vertex set of an equilateral triangle in~$\bf C$,
i.e., is of the form $\alpha\{1,\omega,\omega^2\}$ with $\omega^3=1$.
In~both cases, $\alpha\neq0$ is arbitrary.
\end{defn}

\begin{lem}
\label{lem:main}
In the harmonic case, $L_{\ell,0}$ is a pullback of
$L_{1/2,(2\ell+1)/4,1/4}$, and in the equianharmonic case, $L_{\ell,0}$ is
a pullback of $L_{1/2,1/3,(2\ell+1)/6}$.  Here $\ell\in\bf C$ is arbitrary.
These pullbacks are via maps~$\xi$ which up~to composition with M\"obius
transformations are of the cyclic form $\xi(x)=x^k$, where $k=2,3$
respectively.
\end{lem}

\begin{pf}
The map $\xi_2(x)\defeq x^2$ takes $x=0,\infty$ to $0,\infty$, each with
multiplicity~$2$, and $x=\pm1$ to~$1$ with multiplicity~$1$.  By the theory
of Fuchsian differential operators, any pullback of~$L_{\lambda,\mu,\nu}$
via~$\xi_2$ will be a Fuchsian operator with $-1,0,1,\infty$ as its only
possible singular points.  By Lemma~\ref{lem:exptoexp}, the respective
exponent differences will be $\mu,2\lambda,\mu,2\nu$.  If
$\lambda,\mu,\nu=1/4,1/2,(2\ell+1)/4$, the singular point locations and
exponent differences will be identical to those of~$L_{\ell,B}$ (harmonic
case).  Similarly, any pullback of~$L_{\lambda,\mu,\nu}$ via
$\xi_3(x)\defeq x^3$ will have singular points
$0,1,\omega,\omega^2,\infty$, with exponent differences
$3\lambda,\mu,\mu,\mu,3\nu$.  If $\lambda,\mu,\nu=1/3,1/2,(2\ell+1)/6$, the
point $x=0$ will become an ordinary point, and the singular point locations
and exponent differences will be identical to those of~$L_{\ell,B}$
(equianharmonic case).

The value for the accessory parameter~$B$ of the pullback can be shown to
be zero in both cases.  This follows from Lemma~\ref{lem:firstlem}, since
in both cases a computation (omitted here) yields equal values for the
left-hand and right-hand sides of~(\ref{eq:inlemma}), irrespective
of~$\ell$, iff $B$~is set equal to zero.  It~also follows from a theorem
of~\cite{Maier03}, which determines the values of the accessory parameter
and exponent parameters for which Heun operators are strong pullbacks
of~$L_{\lambda,\mu,\nu}$.

The permutation of $1/4,1/2,(2\ell+1)/4$ into $1/2,(2\ell+1)/4,1/4$, as
required by the statement of the lemma, is accomplished by choosing
$\xi=M\circ\xi_2$, where $M(z)=(z-1)/z$ is the M\"obius transformation that
maps $0,1,\infty$ to $\infty,0,1$.  So in the harmonic case,
$\xi(x)=(x^2-1)/x^2$.  Similarly, the permutation of
$1/3,1/2,{(2\ell+1)/6}$ into $1/2,1/3,{(2\ell+1)/6}$ is accomplished by
composing $\xi_3$ with the map $z\mapsto1-z$.  So in the equianharmonic
case, $\xi(x)=1-x^3$. \qed
\end{pf}

It should be noted that cyclic pullbacks of hypergeometric operators have
been studied or applied by several other authors.  In the harmonic case,
Ivanov \cite{Ivanov2001} discovered that the Jacobi form of the Lam\'e
equation can be reduced to the hypergeometric equation, via a quadratic
transformation analogous to $\xi(x)=x^2$.  In the equianharmonic case,
Clarkson and Olver~\cite{Clarkson96} discovered that the Weierstrass form
of the Lam\'e equation can be similarly reduced, via a cubic transformation
analogous to $\xi(x)=x^3$.  Our efforts to understand their results led
to~\cite{Maier03}, and ultimately to this paper.  Recently, the
Clarkson--Olver transformation has been applied by Kantowski and
Thomas~\cite[Eq.~12]{Kantowski2001}.

\begin{prop}
\label{prop:main}
Let $n$ denote an integer.
\begin{enumerate}
\item In the harmonic case {\rm(}$J=1${\rm)}, $G(L_{n\pm1/6,0})$ is
octahedral if $n\equiv0\ ({\rm mod}~2)$, resp.\ $n\equiv1\ ({\rm mod}~2)$.
\item In the equianharmonic case {\rm(}$J=0${\rm)},
\begin{enumerate}
\item $G(L_{n\pm1/4,0})$ is octahedral if $n\equiv0\ ({\rm mod}~3)$, resp.\
$n\equiv2\ ({\rm mod}~3)$.
\item $G(L_{n\pm1/10,0})$ is icosahedral if $n\equiv0\ ({\rm mod}~3)$, resp.\
$n\equiv2\ ({\rm mod}~3)$.
\item $G(L_{n\pm3/10,0})$ is icosahedral if $n\equiv1\ ({\rm mod}~3)$.
\end{enumerate}
\end{enumerate}
\end{prop}

\begin{pf}
This follows from Lemma~\ref{lem:main}, together with Schwarz's classical
characterization of the $\lambda,\mu,\nu$ for which
$G(L_{\lambda,\mu,\nu})$ is finite.  If the unordered triple
$\lambda,\mu,\nu$ appears on the full Schwarz list, then
$G(L_{\lambda,\mu,\nu})$ will be finite, and the same will be true if a
normalized version of $\lambda,\mu,\nu$ appears there.  Normalization is
performed by replacing $\lambda,\mu,\nu$ by $a\pm\lambda,b\pm\mu,c\pm\nu$,
where $a,b,c$ are any integers whose sum is even.  (See
\cite[Sec.~28]{Poole36}, \cite[Sec.~2.7.2]{Erdelyi53}.)  Preservation of
algebraicity can be verified from the Gauss contiguity relations, which
solutions of hypergeometric equations must satisfy.

The full list includes Cases I, II, IV, VI of Table~\ref{tab:list}, and
also, among others, the icosahedral Case~XIV, for which
$\lambda,\mu,\nu=1/2,2/5,1/3$.  (See \cite[Sec.~30]{Poole36},
\cite[Sec.~2.7.2]{Erdelyi53}.)  By choosing appropriate integers $a,b,c$
(and interchanging the $\mu,\nu$ of Case~XIV), it~is readily verified that
$G(L_{\lambda,\mu,\nu})$ is isomorphic to
\begin{enumerate}
\item $S_4$ if $\lambda,\mu,\nu=1/2,k,1/4$ with $k\in{\bf Z}\pm1/3$ [Case~IV].
\item
\begin{enumerate}
\item $S_4$ if $\lambda,\mu,\nu=1/2,1/3,k$ with $k\in{\bf Z}\pm1/4$ [Case~IV].
\item $A_5$ if $\lambda,\mu,\nu=1/2,1/3,k$ with $k\in{\bf Z}\pm1/5$ [Case~VI].
\item $A_5$ if $\lambda,\mu,\nu=1/2,1/3,k$ with $k\in{\bf Z}\pm2/5$ [Case~XIV].
\end{enumerate}
\end{enumerate}
By Lemma~\ref{lem:main}, each of these $L_{\lambda,\mu,\nu}$ can be pulled
back to a Lam\'e operator of the form $L_{\ell,0}$, with $\ell$ determined
by $(2\ell+1)/4=k$ (Case~1), or by $(2\ell+1)/6=k$ (Cases 2(a)--2(c)).  The
operators~$L_{\ell,0}$ of the proposition are a proper subset: the ones for
which $2\ell\notin\bf Z$.  The reason for imposing this additional
restriction is that if~$2\ell\notin\bf Z$, $G(L_{\ell,0})$ is guaranteed to
be isomorphic to~$G(L_{\lambda,\mu,\nu})$, rather than to a proper
subgroup.  That is because, by Theorem~\ref{thm:main}, the only possible
groups are $S_4$ and~$A_5$, and neither is a subgroup of the other. \qed
\end{pf}

Case~1 of Proposition~\ref{prop:main} provides a counterexample to the
necessary condition of~\cite{Baldassarri81}.  It~should be mentioned that
Case~2(b) is actually a generalization of another result
of~\cite{Baldassarri81}, which is that in the equianharmonic case,
$L_{1/10,0}$ can be pulled back from~$L_{1/2,1/3,1/4}$ via a degree-$3$
cyclic map.  In~fact, Baldassarri was the first to see the relevance of
degree-$3$ cyclic maps in this context.

The following proposition shows that the remaining alternative of
Theorem~\ref{thm:main}, which Proposition~\ref{prop:main} did not cover,
can also be realized.  Unlike Proposition~\ref{prop:main}, it~is specific
to a single value of~$\ell$, and also to a nonzero value of the accessory
parameter~$B$.

\begin{prop}
\label{prop:partial}
Suppose that $J=-80$, i.e., $g_2=80\alpha^2/3$ and $g_3=-80\alpha^3/3$ for
some $\alpha\neq0$; equivalently, that $e_1,e_2,e_3$ are the roots of
$3x^3-20x+20$, multiplied by some $\alpha\neq0$.  Then
$G(L_{1/6,-\alpha/9})$ is icosahedral.
\end{prop}

\begin{pf}
This $\ell=1/6$ example was constructed by a technique suggested by the
method employed to prove Lemma~\ref{lem:main}, which we sketch.  The first
step is to find a rational $\xi(x)$, unramified over ${\bf P}^1({\bf
C})\setminus\{0,1,\infty\}$, such that the pullback of $L_{1/2,1/3,1/5}$
has the same exponent differences as any~$L_{1/6,B}$.  The three singular
points of the pullback that have exponent difference~$1/2$ are taken to be
$e_1,e_2,e_3$.  The second step is to use the formula~(\ref{eq:inlemma}) of
Lemma~\ref{lem:firstlem} to compute the unique~$B$ for which $L_{1/6,B}$
with this choice of $e_1,e_2,e_3$ is, in~fact, the pullback.

It was noted in the proof of Theorem~\ref{thm:main} that in the
$\ell\in{\bf Z}\pm1/6$ icosahedral case, $\xi$~must map the singular point
$x=\infty$ to $z=1$.  Since $x=\infty$ has exponent difference
$\pm(\ell+1/2)=\pm2/3$, the mapping must have multiplicity~2.  In the same
way, it follows that $\xi$~must map each~$e_i$ to $z=0$ with
multiplicity~1.  The function~$\xi$ is characterized by the points in
$\xi^{-1}(\{0,1,\infty\})$ and the multiplicities with which they are
mapped.  Suppose that $\xi^{-1}(0)$ includes $n_0$~ordinary points, besides
$e_1,e_2,e_3$; that $\xi^{-1}(1)$ includes $n_1$~ordinary points,
besides~$\infty$; and that $\xi^{-1}(\infty)$ includes $n_\infty$~ordinary
points.  By Lemma~\ref{lem:exptoexp}, $\xi$~must map each of the
$n_0,n_1,n_\infty$ ordinary points with multiplicity $2,3,5$, respectively.
The integers $n_0,n_1,n_\infty\ge0$ must satisfy
\begin{eqnarray}
\label{eq:firstcond}
&&3+2n_0 = 2 + 3n_1 = 5n_\infty = \deg\xi,\\
&&(3+n_0) + (1+n_1) + n_\infty = 2+\deg\xi.
\label{eq:secondcond}
\end{eqnarray}
Here (\ref{eq:firstcond}) is the degree condition.
Equation~(\ref{eq:secondcond}) is a consequence of the Hurwitz formula,
according to which any rational map $\xi:C\to{\bf P}^1({\bf C})$ from a
nonsingular algebraic curve~$C$ of genus~$g$ to~${\bf P}^1({\bf C})$ that
is unramified above ${\bf P}^1({\bf C})\setminus\{P_1,\ldots,P_r\}$
satisfies $\left|\xi^{-1}(\{P_1,\ldots,P_r\})\right|=2-2g+(r-2)\deg\xi$.

The only solution of (\ref{eq:firstcond})--(\ref{eq:secondcond}) is
$n_0=n_1=n_\infty=1$, with $\deg\xi=5$.  So any function~$\xi$ by which an
operator of the form $L_{1/6,B}$ can be pulled back from $L_{1/2,1/3,1/5}$
must be of the form
\begin{equation}
\label{eq:twopart}
\xi(x) = 
\frac{(x-C_1)(x-C_2)(x-C_3)(x-C_4)^2}{(x-C_5)^5}
=1- \frac{C_6(x-C_7)^3}{(x-C_5)^5}
\end{equation}
for certain $C_1,\ldots,C_7\in\bf C$, where $C_1,C_2,C_3$ are to be
identified with $e_1,e_2,e_3$.  Solutions of~(\ref{eq:twopart}) may be
constructed by elimination theory.  Imposing the condition $e_1+e_2+e_3=0$
yields an essentially unique solution, namely
\begin{equation}
\label{eq:unique}
\xi(x) = \frac{(3x^3-20x+20)(2x-5)^2}{12(x-1)^5}
=1-\frac{(5x-8)^3}{12(x-1)^5},
\end{equation}
which requires $g_2=80/3$ and $g_3=-80/3$.  On the right-hand side
of~(\ref{eq:unique}), $x$~may be replaced by $x/\alpha$ for any
$\alpha\in{\bf C}\setminus\{0\}$.  It~follows by substituting
(\ref{eq:unique}) into~(\ref{eq:inlemma}), and some algebraic manipulation,
that $L_{1/6,B}$ will be a pullback iff $B=-\alpha/9$. The
$\alpha$-dependence is due~to $B$~not being scale-invariant.  \qed
\end{pf}

\begin{cor}
In the nonclassical case $2\ell\notin\bf Z$, finite projective monodromy of
the Lam\'e equation is not uniquely determined by~$\ell$.
\end{cor}
\begin{pf}
By Propositions \ref{prop:main} and~\ref{prop:partial}, $G(L_{1/6,B})$ is
octahedral when $J=1$ and icosahedral when $J=-80$, if in each case, $B$~is
appropriately chosen. \qed
\end{pf}

\section{Explicit Formulas}
\label{sec:explicit}

In practical applications of the Lam\'e equation, such as the astrophysical
application of~\cite{Kantowski2001}, it~is useful to have explicit formulas
for the algebraic solutions, if~any.  The five cases of the following
proposition, which correspond to the four cases of
Proposition~\ref{prop:main} and to Proposition~\ref{prop:partial}, should
serve as examples.

\begin{prop}
\label{prop:examples}
Let $\tau=\tau(x)$, an algebraic complex-valued function of a complex
argument, be defined as follows.
\begin{enumerate}
\item In the harmonic case $\{e_1,e_2,e_3\}=\{-1,0,1\}$, if $\ell=1/6$
and~$B=0$, let $\tau$ be defined by
\begin{displaymath}
\frac{-(\tau^{12}-33\tau^8-33\tau^4+1)^2}{108\,\tau^4(\tau^4-1)^4} =
\frac{x^2-1}{x^2}.
\end{displaymath}
\item In the equianharmonic case $\{e_1,e_2,e_3\}=\{1,\omega,\omega^2\}$,
\begin{enumerate}
\item if $\ell=1/4$ and~$B=0$, let $\tau$ be defined by
\begin{displaymath}
\frac{-(\tau^{12}-33\tau^8-33\tau^4+1)^2}{108\,\tau^4(\tau^4-1)^4} =
1-x^3.
\end{displaymath}
\item if $\ell=1/10$ and~$B=0$, let $\tau$ be defined by
\begin{displaymath}
\frac{[\tau^{30}+522(\tau^{25}-\tau^5)-10005(\tau^{20}+\tau^{10})+1]^2}{1728\,\tau^5(\tau^{10}+11\tau^5-1)^5} =
1-x^3.
\end{displaymath}
\item if $\ell=7/10$ and~$B=0$, let $\tau$ be defined by
\begin{eqnarray*}
&&\frac{[\tau^{30}+522(\tau^{25}-\tau^5)-10005(\tau^{20}+\tau^{10})+1]^2}{1728\,\tau^5(\tau^{10}+11\tau^5-1)^5} \\
&&\qquad\qquad=\frac{s\,(157464\,s^3-352107\,s^2+708750\,s-546875)^2}{(189\,s-125)^5},
\end{eqnarray*}
where $s$ signifies $1-x^3$.
\end{enumerate}
\item In the case when $e_1,e_2,e_3$ are the roots of $3x^3-20x+20$, if
$\ell=1/6$ and $B=-1/9$, let $\tau$ be defined by
\begin{eqnarray*}
&&\frac{[\tau^{30}+522(\tau^{25}-\tau^5)-10005(\tau^{20}+\tau^{10})+1]^2}{1728\,\tau^5(\tau^{10}+11\tau^5-1)^5} \\
&&\qquad\qquad=\frac{(3x^3-20x+20)(2x-5)^2}{12(x-1)^5}.
\end{eqnarray*}
\end{enumerate}
In each of these five cases, the Lam\'e equation $L_{\ell,B}u=0$ has a full
set of algebraic solutions.  Its solution space is spanned by
\begin{equation}
\label{eq:solspace}
\left[\prod_{i=1}^3 (x-e_i)\right]^{-1/4}
\left\{
\frac{1}{\sqrt{d\tau/dx}}
\,,\quad
\frac{\tau}{\sqrt{d\tau/dx}}
\right\},
\end{equation}
where $\tau$ is case-specific.  In cases 1~and~2(a), the projective
monodromy group $G(L_{\ell,B})$, i.e., the Galois group of $\tau$
over~${\bf C}(x)$, is octahedral, and in cases 2(b), 2(c), and~3, it~is
icosahedral.
\end{prop}

\begin{pf}
The solution space~(\ref{eq:solspace}) is of the form specified by
Lemma~\ref{lem:firstlem} in~(\ref{eq:inlemma2}).  In each case, $\tau$~is
defined so that $\tau=\tau'\circ\xi$, where $\xi$~is the rational function
by which $L_{\ell,B}$ is pulled back from some~$L_{\lambda,\mu,\nu}$, and
$\tau'$~is a ratio of solutions of~$L_{\lambda,\mu,\nu}v=0$.  In~all cases
except~2(c), the right-hand side of the defining equation is~$z=\xi(x)$, as
supplied in the proof of Lemma~\ref{lem:main} or the proof of
Proposition~\ref{prop:partial}, and the left-hand side is the appropriate
polyhedral function, as supplied in the final column of
Table~\ref{tab:list}, applied to~$\tau$.

Case~2(c) is special.  As was sketched in the proof of
Proposition~\ref{prop:main}, $L_{7/10,0}$ is the pullback via
$\xi(x)=1-x^3$ of $L_{1/2,1/3,2/5}$, which is Schwarz's Case~XIV (modulo
the interchange of $\mu,\nu$).  Case~XIV is not on the basic Schwarz list,
and in~fact, it is not the case that a ratio $\tau'=\tau'(z)$ of
independent solutions of $L_{1/2,1/3,2/5}v=0$ is the inverse of a rational
function.  However, Case~XIV is itself a pullback of the basic icosahedral
Case~VI.  So one can choose $\tau'=\tau''\circ\xi'$, where $\xi'$~is the
rational function by which $L_{1/2,1/3,2/5}$ is pulled back from
$L_{1/2,1/3,1/5}$, and $\tau''$~is a ratio of solutions of
$L_{1/2,1/3,1/5}v=0$, the inverse of which is listed in
Table~\ref{tab:list}.  The formula in Case~2(c) defines $\tau$ so that
$\tau=\tau'\circ\xi=\tau''\circ\xi'\circ\xi$.

A rational map~$\bar\xi$ equivalent to~$\xi'$ was worked~out by Klein in
1877, in a paper in which he completed the reduction of the Schwarz list to
the basic Schwarz list~\cite[Sec.~10]{Klein1877}.  His formula was
\begin{equation}
\bar\xi(s) = 1-
\frac{(64\,s+189)(64\,s^2+133\,s + 49)^3}
{7^7\cdot27\cdot(s+1)^2}
\end{equation}
which maps $s=0,-189/64,-1$ respectively to~$\bar\xi(s)=0,1,\infty$.  For
our purposes, this morphism must be composed with a M\"obius
transformation.  Composing with $M(s)=189s/(125-189s)$, which takes
$s=0,1,\infty$ to $M(s)=0,-189/64,-1$, yields
\begin{equation}
(\bar\xi\circ M)(s) = 
\frac{s\,(157464\,s^3-352107\,s^2+708750\,s-546875)^2}{(189\,s-125)^5}
\end{equation}
as the rational map~$\xi'=\xi'(s)$ by which $L_{1/2,1/3,2/5}$ is pulled
back from $L_{1/2,1/3,1/5}$.  This map appears on the right-hand side in
Case~2(c).  \qed
\end{pf}

\section{The Weierstrass Form}
\label{sec:corrections}

In classical treatments~\cite{Whittaker27}, the Weierstrass-form Lam\'e
equation is regarded as an equation on~$\bf C$, of the form
\begin{equation}
\label{eq:old}
\frac{d^2u}{dt^2} - \left[\ell(\ell+1)\wp(t)+B\right]u = 0.
\end{equation}
Here $\wp:{\bf C}\to{\bf P}^1({\bf C})$ is the Weierstrass $\wp$-function
corresponding to some period lattice ${\mathcal L}={\omega_1}{\bf Z}+
{\omega_2}{\bf Z}\subset{\bf C}$, with $\omega_1,\omega_2$ independent
over~$\bf R$; i.e., $(\wp')^2 = 4\wp^3-g_2\wp-g_3$ for some $g_2,g_3\in\bf
C$ for which $\Delta=g_2^3-27g_3^2\neq0$.  Equation~(\ref{eq:old}) is a
Schr\"odinger equation with an elliptic potential, extended to the complex
domain.  The algebraic Lam\'e equation $L_{\ell,B}u=0$ can be obtained
from~(\ref{eq:old}) by the substitution $x=\wp(t)$.  That~is,
(\ref{eq:old})~is the strong pullback to~$\bf C$ of $L_{\ell,B}u=0$
by~$\wp$.

Another interpretation is possible.  The map $\wp:{\bf C}\to{\bf P}^1({\bf
C})$ is the composition of two maps, $\phi:{\bf C}\to E_{g_2,g_3}$ and
$\pi:E_{g_2,g_3}\to{\bf P}^1({\bf C})$.  Here $E_{g_2,g_3}$ is the elliptic
curve specified by $y^2=4x^3-g_2x-g_3$, and the maps $\phi$ and~$\pi$ are
defined by $\phi(t)=\left(\wp(t),\wp'(t)\right)$ and $\pi(x,y)=x$.
$E_{g_2,g_3}$ is homeomorphic to a torus, and the projection~$\pi$ is a
double cover of ${\bf P}^1({\bf C})$ by~$E_{g_2,g_3}$.  From an
algebraic-geometric point of view, it~is more reasonable to pull the
algebraic-form Lam\'e equation back to~$E_{g_2,g_3}$ via~$\pi$, than
to~$\bf C$ via $\pi\circ\phi$.  We~call the resulting equation
on~$E_{g_2,g_3}$ the Weierstrass-form Lam\'e equation, and write it
$L_{\ell,B,g_2,g_3}u=0$.  By~examination, the operator $L_{\ell,B,g_2,g_3}$
has only one singular point, namely the point~$O$, i.e.,
$(x,y)=(\infty,\infty)$, where its characteristic exponents are
$-\ell,\ell+1$.  We shall informally regard $E_{g_2,g_3}$ as a subset of
${\bf P}^1({\bf C})\times {\bf P}^1({\bf C})$, coordinatized by~$(x,y)$,
although in a more careful treatment $E_{g_2,g_3}$~would be defined as the
projective curve $y^2z=4x^3-g_2xz^2-g_3z^3$ in~${\bf P}^2({\bf C})$,
equipped with homogeneous coordinates~$(x,y,z)$.

The pullback theory of Section~\ref{sec:prelims} applies when the algebraic
curve~$C$ equals $E_{g_2,g_3}$, just as it applied when $C={\bf P}^1({\bf
C})$ and $K={\bf C}(x)$.  The function field~$\tilde K$ on~$E_{g_2,g_3}$ is
${\bf C}(x,y)\supset{\bf C}(x)$, a~degree-2 extension, and the derivation
$D=d/dx$ extends in the obvious way to~$\tilde K$, via
$Dy\defeq(12x^2-g_2)/2y$.  With these choices, $L_{\ell,B,g_2,g_3}$ is of
the form $D^2 + \tilde{\mathcal A}\cdot D + \tilde{\mathcal B}$, for
$\tilde{\mathcal A},\tilde{\mathcal B}\in K\subset\tilde K$.  Since the
Wronskian is algebraic, the projective monodromy group
$G(L_{\ell,B,g_2,g_3})$ is finite iff a ratio of solutions~$\tilde\tau$ of
$L_{\ell,B,g_2,g_3}u=0$ on~$E_{g_2,g_3}$ is algebraic over ${\bf C}(x,y)$.
This is equivalent to $G(L_{\ell, B})$ being finite, which occurs iff a
ratio of solutions~$\tau$ of $L_{\ell,B}u=0$ on~${\bf P}^1({\bf C})$ is
algebraic over~${\bf C}(x)$.  The equivalence is due~to the (local)
pullback property $\tilde\tau=\tau\circ\pi$, which implies that the two
sorts of algebraicity are equivalent.  But the two groups may not be
isomorphic.  This is because $\tilde\tau$~may have lower degree over ${\bf
C}(x,y)$ than $\tau$~has over~${\bf C}(x)$.

In a thought-provoking paper, Churchill~\cite{Churchill99} examined the
implications of finite group theory for the monodromy
of~$L_{\ell,B,g_2,g_3}$, as~well as for the monodromy of the hypergeometric
operator $L_{\lambda,\mu,\nu}$.  The (projective) monodromy group
of~$L_{\lambda,\mu,\nu}$ is doubly generated: it~is generated by the images
of loops around any two of the singular points $z=0,1,\infty$.  Similarly,
since $E_{g_2,g_3}$ is homeomorphic to a torus, the (projective) monodromy
group of $L_{\ell,B,g_2,g_3}$ is generated by the images of only two loops.
But the conjugacy classes of the two monodromy (resp.\ projective
monodromy) generators in ${\it GL}(2,{\bf C})$ (resp.\ ${\it PGL}(2,{\bf
C})$) are determined by the characteristic exponents of the singular
point(s).  This constrains what, up~to isomorphism, the monodromy group
(resp.\ projective monodromy group) may be, in the case when it is finite.

In this way, Churchill was able to obtain significant results on the
projective monodromy of~$L_{\ell,B,g_2,g_3}$ without using pullbacks.
Like~$G(L_{\ell,B})$, $G(L_{\ell,B,g_2,g_3})$ can never be cyclic, and can
be dihedral only if~$2\ell\in\bf Z$.  Moreover, in the nonclassical case
$2\ell\notin\bf Z$, it~cannot be dihedral.  He~showed that in the
nonclassical case, $G(L_{\ell,B,g_2,g_3})$ can be tetrahedral only if
$\ell\in{\bf Z}\pm1/4$, and can be octahedral or icosahedral only if
$\ell\in{\bf Z}\pm1/10$, $\ell\in{\bf Z}\pm1/6$, or $\ell\in{\bf
Z}\pm3/10$.

Necessarily $G(L_{\ell,B,g_2,g_3})\trianglelefteq
G(L_{\ell,B})$~\cite{Churchill99}, so conditions on~$G(L_{\ell,B,g_2,g_3})$
yield conditions on~$G(L_{\ell,B})$, and vice versa.  By combining his
results with those of~\cite{Baldassarri81}, Churchill was able to deduce
that $G(L_{\ell,B,g_2,g_3})$ cannot be octahedral.  Unfortunately this is
incorrect: the proofs of Theorem~5.3 and Corollary~5.4
of~\cite{Churchill99}, which include this assertion, rely crucially on the
incorrect result of~\cite{Baldassarri81} that $G(L_{\ell,B})$ can be
octahedral only if $\ell\in{\bf Z}\pm1/4$.

Theorem~\ref{thm:mainmain} is a characterization of $G(L_{\ell,B,g_2,g_3})$
as~well as $G(L_{\ell,B})$, which is obtained from pullback theory alone.
This theorem builds~on and subsumes Theorem~\ref{thm:main}.

\begin{thm}
\label{thm:mainmain}
The equation $L_{\ell,B,g_2,g_3}u=0$ on~$E_{g_2,g_3}$ has a full set of
algebraic solutions iff $G(L_{\ell,B,g_2,g_3})$ is finite, which is
equivalent to $L_{\ell,B}u=0$ on~${\bf P}^1({\bf C})$ having a full set of
algebraic solutions, and to $G(L_{\ell,B})$ being finite.  In the
nonclassical case $2\ell\notin\bf Z$, the following are the only ways this
can occur.
\begin{enumerate}
\item $G(L_{\ell,B,g_2,g_3})$ is tetrahedral and $G(L_{\ell,B})$ is
octahedral; in which case $\ell$~must equal $n\pm1/4$, with $n$~an integer.
\item $G(L_{\ell,B,g_2,g_3})$ is octahedral and $G(L_{\ell,B})$ is
octahedral; in which case $\ell$~must equal $n\pm1/6$, with $n$~an integer.
\item $G(L_{\ell,B,g_2,g_3})$ is icosahedral and $G(L_{\ell,B})$ is
icosahedral; in which case $\ell$~must equal $n\pm1/10$, $n\pm1/6$, or
$n\pm3/10$, with $n$~an integer.
\end{enumerate}
All five of the preceding alternatives can be realized.
\end{thm}

\begin{pf}
The first sentence has already been proved.  The proof of the necessary
conditions on~$\ell$ for $G(L_{\ell,B,g_2,g_3})$ to be each possible finite
group is similar to the proof of Theorem~\ref{thm:main}: it~runs down the
rows of the basic Schwarz list, beginning with the tetrahedral.  Before
beginning the proof, note that the pullback function~$\xi$ guaranteed to
exist by Theorem~\ref{thm:pullback} will map the lone singular point~$O$
of~$L_{\ell,B,g_2,g_3}$ to one of $\{0,1,\infty\}$, since the exponent
difference at~$O$ is $\pm(2\ell+1)$, and an analogue of
Lemma~\ref{lem:lookahead} holds.  Also, note that the degree
formula~(\ref{eq:degreeformula}) of Baldassarri and Dwork yields
$\pm\ell=(\deg\xi)/3$ when applied to $F=L_{\ell,B,g_2,g_3}$ and
$F'=L_{1/2,1/3,1/3}$, and $\pm\ell=(\deg\xi)/6$ when applied to
$F=L_{\ell,B,g_2,g_3}$ and $F'=L_{1/2,1/3,1/4}$, since $E_{g_2,g_3}$, being
elliptic, has genus~$g=1$.

If $G(L_{\ell,B,g_2,g_3})$ is tetrahedral, $L_{\ell,B,g_2,g_3}$ must be a
pullback of $L_{1/2,1/3,1/3}$.  By Lemma \ref{lem:exptoexp}, if $\xi(O)=0$
then $2\ell+1$ is an integer multiple of~$1/2$, i.e., $\ell\in{\bf
Z}\pm1/4$.  The possibilities $\xi(O)=1,\infty$ can be ruled~out, since
they would imply respectively that $\xi^{-1}(\infty),\xi^{-1}(1)$ consists
of ordinary points, each mapped with multiplicity~3.  Either would imply
$3\mid\deg\xi$, which with $\pm\ell=(\deg\xi)/3$ would contradict
$2\ell\notin\bf Z$.

If $G(L_{\ell,B,g_2,g_3})$ is octahedral, $L_{\ell,B,g_2,g_3}$ must be a
pullback of $L_{1/2,1/3,1/4}$.  By Lemma \ref{lem:exptoexp}, if $\xi(O)=1$
then $2\ell+1$ is an integer multiple of~$1/3$, i.e., $\ell\in{\bf
Z}\pm1/6$.  The possibilities $\xi(O)=0,\infty$ can be ruled~out.  If
$\xi(O)=0$ then $\xi^{-1}(1),\xi^{-1}(\infty)$ consist of ordinary points,
each mapped with multiplicity~3,4, respectively.  This would imply
$3\mid\deg\xi$ and $4\mid\deg\xi$, hence $12\mid\deg\xi$; which with
$\pm\ell=(\deg\xi)/6$ would contradict $2\ell\notin\bf Z$.  $\xi(O)=\infty$
is ruled~out similarly.

If $G(L_{\ell,B,g_2,g_3})$ is icosahedral, $L_{\ell,B,g_2,g_3}$ must be a
pullback of $L_{1/2,1/3,1/5}$.  By Lemma \ref{lem:exptoexp}, if $\xi(O)=1$
then $2\ell+1$ is an integer multiple of~$1/3$, i.e., $\ell\in{\bf
Z}\pm1/6$, and if $\xi(O)=\infty$ then $2\ell+1$ is an integer multiple
of~$1/5$, i.e., $\ell\in{\bf Z}\pm1/10$ or $\ell\in{\bf Z}\pm3/10$.  The
possibility $\xi(O)=0$ can be ruled~out, since it would imply that
$2\ell+1$ is an integer multiple of~$1/2$, i.e., $\ell\in{\bf Z}\pm1/4$.
The group $G(L_{\ell,B})$ is finite if $G_{\ell,B,g_2,g_3}$ is finite, so
if $\xi(O)=0$, Theorem~\ref{thm:main} implies that $G(L_{\ell,B})$ is
octahedral.  But $G(L_{\ell,B,g_2,g_3})$ must be isomorphic to a subgroup
of $G(L_{\ell,B})$.

The classification scheme of the theorem results from combining the
just-derived conditions on~$G(L_{\ell,B,g_2,g_3})$ with the conditions of
Theorem~\ref{thm:main} on $G(L_{\ell,B})$.  That $G(L_{\ell,B,g_2,g_3})$
octahedral implies $G(L_{\ell,B})$ octahedral is due~to $A_5$ not having
any normal $S_4$ subgroup.  The realizability of all five alternatives was
proved in Section~\ref{sec:result} (it~follows from Propositions
\ref{prop:main} and~\ref{prop:partial}). \qed
\end{pf}

By Theorem~\ref{thm:mainmain}, $G(L_{\ell,B})$ being octahedral does not
uniquely determine the group $G(L_{\ell,B,g_2,g_3})$: it~may be either
octahedral or tetrahedral.  The latter occurs when the extension ${\bf
C}(x,y,\tilde\tau)/{\bf C}(x,y)$ has lower degree than ${\bf
C}(x,\tau)/{\bf C}(x)$.  The two possibilities are exemplified by Cases
1~and~2(a) of Proposition~\ref{prop:examples}, respectively, which have
$\ell,B,g_2,g_3$ equal to $1/6,0,4,0$ and $1/4,0,0,4$.  A~ratio
$\tilde\tau$ of solutions of $L_{\ell,B,g_2,g_3}u=0$ is specified by
\begin{equation}
\label{eq:cases}
\frac{-(\tilde\tau^{12}-33\tilde\tau^8-33\tilde\tau^4+1)^2}{108\,\tilde\tau^4(\tilde\tau^4-1)^4} =
\cases
{
{\displaystyle\frac{x^2-1}{x^2}}, & $\ell,B,g_2,g_3 = 1/6,0,4,0$;\cr
1-x^3, & $\ell,B,g_2,g_3 = 1/4,0,0,4$,\cr
}
\end{equation}
since $\tilde\tau=\tau\circ\pi$.  In the second case,
$y^2=4x^3-g_2x-g_3=4x^3-4$, so ${1-x^3}=-y^2/4$, implying that the minimum
polynomial of $\tau$ over~${\bf C}(x)$ is reducible over~${\bf C}(x,y)$.
In~fact, $\tilde\tau$ can be chosen to satisfy
\begin{equation}
\frac{\tilde\tau^{12}-33\tilde\tau^8-33\tilde\tau^4+1}{(\sqrt{108}/2)
\tilde\tau^2(\tilde\tau^4-1)^2} = \pm y,
\end{equation}
where either sign is acceptable. Each sign yields a 12-branched algebraic
function~$\tilde\tau$ on the equianharmonic elliptic curve~$E_{0,4}$ (with
$J=0$) that projects to~$\tau$, rather than yielding a 24-branched function
on~$E_{0,4}$.  And ${\bf C}(x,y,\tilde\tau)$ is obtained from~${\bf C}(x)$
via the tower ${\bf C}(x)\subset{\bf C}(x,y)\subset{\bf
C}(x,y,\tilde\tau)$, where the extensions are algebraic of degrees 2
and~12, respectively.  The group $G(L_{1/4,0,0,4})$, which is the Galois
group of $\tilde\tau$ over~${\bf C}(x,y)$, has order~12 and must be
tetrahedral, i.e., isomorphic to~$A_4$.

In general, this reduction may not occur.  In the first case
of~(\ref{eq:cases}), in which $y^2=4x^3-4x$, the analogous substitution
does not lead to a reduction of the degree.  The function~$\tilde\tau$ on
the harmonic elliptic curve~$E_{4,0}$ (with $J=1$) is 24-branched,
like~$\tau$, the function to which it projects.  So the group
$G(L_{1/6,0,4,0})$ has order~24 and must be octahedral, i.e., isomorphic
to~$S_4$.

\ack 
The author gratefully acknowledges the hospitality of the Texas
Institute for Computational and Applied Mathematics (TICAM).


\small\def\em{\it} \newcommand{\noopsort}[1]{} \newcommand{\printfirst}[2]{#1}
  \newcommand{\singleletter}[1]{#1} \newcommand{\switchargs}[2]{#2#1}

\end{document}